\documentclass{ip-journal}
\usepackage[utf8]{inputenc}
\usepackage{amsmath,amssymb,graphicx, xypic,amsthm}
\title{(0,2)-Deformations and the $G$-Hilbert Scheme}
\author{Benjamin Gaines}

\newcommand\CC{\mathbb{C}}

\newcommand\cO{\mathcal{O}}

\renewcommand\cD{\mathcal{D}}
\renewcommand\cH{\mathcal{H}}
\newcommand\cT{\mathcal{T}}

\newcommand\cW{\mathcal{W}}

\newcommand\RR{\mathbb{R}}
\newcommand\ZZ{\mathbb{Z}}
\renewcommand\AA{\mathbb{A}}

\makeatletter
\newtheorem*{rep@theorem}{\rep@title}
\newcommand{\newreptheorem}[2]{%
\newenvironment{rep#1}[1]{%
 \def\rep@title{#2 \ref{##1}}%
 \begin{rep@theorem}}%
 {\end{rep@theorem}}}
\makeatother

\newtheorem{theorem}{Theorem}[section]
\newtheorem{lemma}[theorem]{Lemma}
\newtheorem{proposition}[theorem]{Proposition}

\newreptheorem{theorem}{Theorem}
\newreptheorem{lemma}{Lemma}

\DeclareMathOperator\ext{Ext}

\DeclareMathOperator\spec{Spec}

\begin{document}

\begin{abstract}
We study first order deformations of the tangent sheaf of resolutions of Calabi-Yau threefolds that are of the form $\CC^3/\ZZ_r$, focusing
 on the cases where the orbifold has an isolated singularity.  We prove a lower bound on the number 
of deformations for any crepant resolution of this orbifold.  We show that this lower bound is achieved when the resolution used is the 
$G$-Hilbert scheme, and note that this lower bound can be found using methods from string theory.  These methods lead us to a new way to construct the 
$G$-Hilbert scheme using the singlet count.

\end{abstract}
\maketitle

\section{Introduction}

When $G$ is an abelian group, the orbifold $\CC^3/G$ is a normal toric variety.  It can have several crepant resolutions, and the process for finding these 
resolutions is explained in detail in \cite{Toricbook}.  We will review these methods in Section 2.  For each of these resolutions, we are interested in 
the moduli space of the tangent bundle $\cT$.  We will calculate the dimension of $\ext^1(\cT,\cT)$, to determine the number of first order deformations 
of the tangent sheaf.

There are many reasons to be interested in this moduli space, but our primary motivation comes from string theory.  One of the predictions that 
geometric string theory makes is that there is connection between the orbifold conformal theory, independent of the projective resolution chosen, and the 
moduli space.  In $N=(2,2)$ theories, this statement is well understood and follows from the McKay correspondence, as in \cite{BatyrevDais,BKR}.  
In the case of $(0,2)$-theories, this relationship is less well understood.  One of our goals is to look for a relationship between the 
conformal field theory and the number of deformations of the tangent bundle of $X$, as predicted in \cite{DHVW}.  
Following Aspinwall, we call this the $(0,2)-$McKay correspondence \cite{Paul11}.

In this paper, we prove that the conformal field theory orbifold count will serve as a lower bound for the dimension of the moduli space.  We 
also find that the conformal field theory gives a distinguished quiver, and that this quiver can be used to construct a particular resolution for $X$.  
We will then show that this particular resolution is in fact the $G$-Hilbert scheme, and that for this resolution the $(0,2)-$McKay correspondence holds.
  This gives the following main result:
\begin{theorem}
\label{maintheorem}
 Among all crepant resolutions of $\CC^3/\ZZ_r$ with $r$ prime, the $G$-Hilbert scheme has minimal dimension for $\ext^1(\cT,\cT)$.  This dimension exactly 
matches the prediction given by the conformal field theory of the orbifold. 
\end{theorem}

In the course of proving this theorem, we find a few other interesting results about the $G$-Hilbert scheme, in the case where 
$G$ is an abelian subgroup of $SL(3,\CC)$.  In particular, we see that there is a very easy way to find the number of first order deformations of the 
$G$-Hilbert scheme, using the knockout method of \cite{CR}.  We also show that the fan for the $G$-Hilbert scheme can be constructed using the singlets 
from the conformal field theory.

\section{Crepant Resolution of an Orbifold}
\subsection{Quotient Presentation of Toric Variety}
\indent We begin with a review of the quotient presentation of a toric variety.  Let $\Sigma$ be a simplicial fan, 
with 3 dimensional lattice $N$.
Let $\{\rho_i\}$ denote the set of one-dimensional cones (or rays) of $\Sigma$, with their minimal generators 
given by $u_{\rho_i}$, and $|\Sigma(1)|$ being the number of rays.  Let $M$ be the dual lattice to $N$, and 
$\{P_i\}$ be the set of primitive collections of $\Sigma$.  A primitive collection is  a set of rays so that 
the entire set does not lie in a cone of $\Sigma$, but every proper subset does.  
We can now describe the toric variety associated with $\Sigma$ as a quotient $X=(\spec(S)- V(B))/ G$, 
using the method of Cox \cite{Cox}. In this notation, $S$ is the homogeneous coordinate ring of the variety, $B$ is the irrelevant 
ideal, and $G$ is the quotient group.  
\begin{itemize}
 \item $S=\CC[x_{\rho_i}]$, where each ray $\rho_i$ corresponds to a variable $x_{\rho_i}$.  
 \item $G=\text{Hom}_\ZZ(Cl(X),\CC^*)$, with a group action given by 
$$\{(t_{\rho})\in (\CC^*)^{|\Sigma(1)|} |\prod t_{\rho}^{<e_i,u_{\rho}>}=1\text{ for }1\leq i \leq 3\}$$
where $\{e_i\}$ are a basis for  $M$.
 \item $ B=\bigcap(I_{P_i})$, where $I_{P_i}$ is the ideal generated by the elements of the primitive collection $P_i$.
\end{itemize}
In particular, we note that both $S$ and $G$ depend only on the one-dimensional cones of $\Sigma$.

\subsection{Exceptional Divisors}
Let $X=\CC^3/G$, with $G=\ZZ_r$ and having a group action given by $(\zeta_r^{a_1},\zeta_r^{a_2},\zeta_r^{a_3})$ where
$\zeta_r$ is a primitive $r$-th root of unity.  We are interested in resolving $\CC^3/G$, and getting a smooth variety. 
To begin, we'll look at the case where 
$a_1=1$, $a_2=a$, and $a_3=b$.  As long as some $a_i$ and $r$ are relatively prime, we can renumber to get this case.  We note that if $X$ has an 
isolated singularity, this is always true.  
Since $\CC^3/G$ is a simplicial toric variety, there is a fan in $\RR^3$ associated to it.  Such a fan must consist of 
a single cone, with $3$ rays.  To find minimal generators for these rays, we must make sure the group action described 
in the previous section gives $(\zeta_r^{a_1},\zeta_r^{a_2},\zeta_r^{a_3})$.  We do this by choosing 
$u_1=(r,-a,-b)$, $u_2=(0,1,0)$, and $u_3=(0,0,1)$, where the lattice $N$ is just $\ZZ^3$.\\   
\indent We must now find all of the lattice points in the interior of the convex hull of these 3 minimal generators.  We note 
that these points all lie in a hyperplane, and therefore can be used to find a crepant resolution of $X$. 
Each of these points will be the minimal generator of a ray in $\RR^3$, which corresponds to an exceptional divisor of the resolved variety.
\begin{proposition}
 The interior points can be labeled by $$\{(i,-s,-t)| i,s,t\in \ZZ,0<i<r,rs<ia<r(s+1),rt<ib<r(t+1),i-s-t=1\}$$
\end{proposition}
  We will now prove the proposition.  
We can find interior points by starting at one vertex, $(r,-a,-b)$ and adding $d_1(-r,1+a,b)$ and $d_2(-r,a,1+b)$, satisfying 
$d_1+d_2\leq 1$, $d_1,d_2\geq 0$, ignoring the cases corresponding to vertices.
%In other words, we're adding a combination of the rays that go from the starting vertex to the 2 other vertices, which 
%are not strictly along those rays.  
We therefore see that every interior point is of the form $(r(1-d_1-d_2),(d_1+d_2-1)a+d_1,(d_1+d_2-1)b+d_2)$.  To find just the lattice points, we need 
all of these values to live on the lattice $N$.  In this case, the lattice was $\ZZ^3$.  Let $r(1-d_1-d_2)=i$.  
%It is not difficult to show that the set of all interior points is given by 
%$$\{(r(1-c_1-c_2),(c_1+c_2-1)a+c_1,(c_1+c_2-1)b+c_2) : c_1+c_2<1,c_1>0,c_2>0\}$$
%Since we want lattice points, we must have 
%see that $1-(c_1+c_2)=\frac{i}{r}$.  We can plug this into the above equation for interior points to see all interior lattice points are of the form 
%$(i,\frac{-ia}{r}+c_1,\frac{-ib}{r}+c_2)$.  Since $0<c_1<1$, and $\frac{-ia}{r}+c_1$ is an integer, we know that $\frac{-ia}{r}+c_1$ is equal to 
%$-[\frac{ia}{r}]$, where $[n]$ is the smallest integer less than $n$, and 
After a few algebraic manipulations, this leads to the conditions that $d_1r\equiv ia \mod r$, 
%.  Similarly, $\frac{-ib}{r}+c_2=-[\frac{ib}{r}]$ and 
$d_2r\equiv ib \mod r$, and
%  So an interior point 
%must look like $(i,-[\frac{ia}{r}],-[\frac{ib}{r}])$.  However, there's another condition we also must include.  We needed to have 
$d_1+d_2=1-\frac{i}{r}$.  %Using the above, this is the same as saying $\frac{i}{r}+\frac{ai}{r} \mod(1)+ \frac{bi}{r} \mod(1)=1$.
We therefore have $d_1=\frac{ia}{r}-s$, $d_2=\frac{ib}{r}-t$, for some $s,t$ integers.  

Thus, the set of interior lattice points is labeled by $$\{(i,-s,-t)| i,s,t\in \ZZ,0<i<r,rs<ia<r(s+1),rt<ib<r(t+1),i-s-t=1\}$$  We will call these interior points $u_4$ 
through $u_N$ (and will call the three vertices $u_1,u_2,u_3$).$\square$\\ \\

As an example, let's look at the orbifold $X=\CC^3/\ZZ_{11}$, where the group action is given by 
$(\zeta_{11},\zeta_{11}^{2},\zeta_{11}^{8})$.  In this case, the minimal generators for the rays of the cone are 
$u_1=(11,-2,-8), u_2=(0,1,0),$ and $u_3=(0,0,1)$.  Using the process described above, we must find the set 
$$\{(i,-s,-t)|i,s,t\in \ZZ, 11s<2i<11s+11, 11t<8i<11t+11, i-s-t=1\}.$$  Checking all possible values for $i$, we 
find that the interior points are $$u_4=(1,0,0), u_5=(2,0,-1), u_6=(3,0,-2), u_7=(6,-1,-4), u_8=(7,-1,-5)$$
These can be seen in the figure below.
\begin{center}
\includegraphics[height=2in]{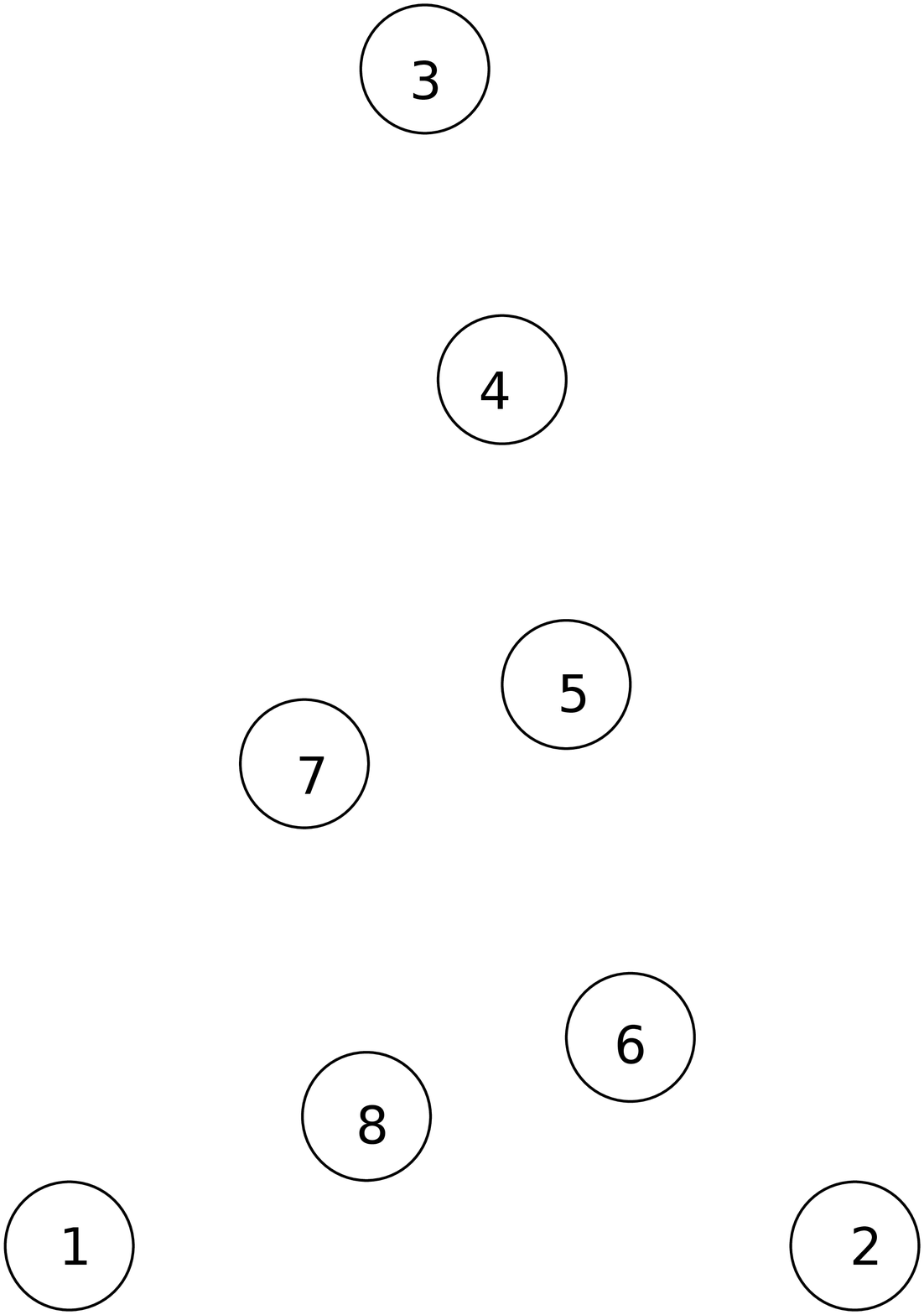}\\
Figure 2.1:  $\CC^3/\ZZ_{11}$ with Interior Lattice Points\\
\end{center}

Note that if the singularity is not isolated, it is possible for the group action to be given by $(\zeta_r^{a_1},\zeta_r^{a_2},\zeta_r^{a_3})$, where 
each of the $a_i$ is not relatively prime to $r$.  In this case a similar theorem holds (by finding interior points on a sublattice of $\ZZ^3$), 
but for reasons we will expound on later, we will limit our focus to orbifolds with isolated singularities.

%\subsection{Non-Isolated Case}
%Now say the group action is given by  $(\zeta_r^{a_1},\zeta_r^{a_2},\zeta_r^{a_3})$ with no condition on $a_1$.  Constructing the 
%cone corresponding to this toric variety is a little trickier. We would like the minimal generators for the rays to be 
%$u_1=(r,-a_2,-a_3),u_2=(0,a_1,0),$ and $u_3=(0,0,a_1)$, in order for $G$ to have the correct group action.  However this means the 
%lattice cannot be $\ZZ^3$, or $u_2,u_3$ would not be minimal generators.  By instead making the lattice $N=(a_1\ZZ)^3$, 
%we can make it so that these are in fact minimal generators for their rays.  The basis 
%basis for lattice: (0,a_1,0),(0,0,a_1),(a_1,1,1)
%Now the fan $\Sigma$ corresponds to a toric variety with the appropriate group action.  
% We now want to find all interior lattice points, as we did previously.  We can use the same method as described above, but 
%rather than require each coordinate be an integer, we require the point to live on the new lattice.  
%This gives a subset of the set of points found using the above method, and this subset of points gives exactly the 
%rays needed to crepantly resolve $X$.

\subsection{Resolving the Singularity}
Now that we have found all of the exceptional divisors, we can find a resolution of the singularity.  To do this, 
we find some triangulation of the points $u_1,..u_N$, and then take the fan over this triangulation.  This gives 
a new fan with the same lattice $N$, which corresponds to the crepant resolution associated with this particular 
triangulation.  Since all of these fans have the same one dimensional rays (with minimal generators $\{u_1,..,u_N\}$), 
they have the same homogeneous coordinate ring $S$ and quotient group $G$.
 We have $S=\CC[x_1,\ldots,x_N]$, where each $x_i$ corresponds to a ray.  To find $G$, 
we first construct an $n\times 3$ matrix $A$, whose rows are the coordinates of these points.
$$A=\left(\begin{array}{c c c}
r & -a & -b \\
0 & 1 & 0\\
0 & 0 & 1\\
(u_4)_1 & (u_4)_2 & (u_4)_3 \\
\vdots & \vdots & \vdots \\
%\nu^\alpha_1r & -\nu^\alpha_1a+\nu^\alpha_2 & -\nu^\alpha_1b+\nu^\alpha_3 \\
(u_N)_1 & (u_N)_2 & (u_N)_3
\end{array}\right)$$ 
The cokernel of $A$ will give the class group $Cl(X)$ of any resolution with these one-dimensional cones, and $G=Hom_\ZZ(Cl(X),\CC^*)$.  This 
tells us that $G\cong (\CC^*)^{N-3}$, and what the group action of $G$ is.  Let $\Phi$ be a matrix with coefficients in 
$\ZZ$ and rows given by a basis of the cokernel of $A$.  In the physics literature, this is known as the charge matrix.  
This matrix will be $(N-3)\times N$, and so using some linear algebra we can reduce it to another matrix $\Phi'$, where 
%singletcohomologyequivalence.tex!!!
$$\Phi'=\left(\begin{array}{c c c c c c c}
\nu_1^4 & \nu_2^4 & \nu_3^4 & -1 & 0 & \cdots 	& 0 \\
\nu_1^5 & \nu_2^5 & \nu_3^5 & 0 & -1 & \cdots 	& 0\\
  &  & &\vdots & &\\
\nu_1^N & \nu_2^N & \nu_3^N & 0 & 0 & \cdots & -1 \\
%\nu^\alpha_1r & -\nu^\alpha_1a+\nu^\alpha_2 & -\nu^\alpha_1b+\nu^\alpha_3 \\
\end{array}\right)$$ 
%so each column is $(\nu_1^i,\nu_2^i,\nu_3^i,0...,0,-1,0...0)$, with $-1$ in the $i$th row, 
for some rational numbers 
$\nu_1^i,\nu_2^i,\nu_3^i$ (note:  the superscript is an index here, not a power).  

The matrix $\Phi$ gives a grading for $S$, where the multidegree for $x_j$ is a 
vector given by the $j$th column of $\Phi$.  We will call this vector the charge of $x_j$, and denote it by $q_j$.\\ 

%In particular, $G$ will act on $x_i$ by sending it to CHARGES!!!!

Let's go back to our example, of $\CC^3/\ZZ_{11}$.  Using the result from the previous section, we know there are 
$8$ rays in the fan of any crepant resolution, so $S=\CC[x_1,..,x_8]$.    We also find that 
$$A=
\left(\begin{array}{c c c}
11 & -2 & -8 \\
0 & 1 & 0\\
0 & 0 & 1\\
1 & 0 & 0 \\
2 & 0 & -1 \\
3 & 0 & -2 \\
6 & -1 & -4 \\
7 & -1 & -5 
\end{array}\right) \qquad
\Phi'=\left(\begin{array}{c c c c c c c c}
\frac{1}{11} & \frac{2}{11} & \frac{8}{11} & -1 & 0 & 0 & 0 & 0  \\
\frac{2}{11} & \frac{4}{11} & \frac{5}{11} & 0 & -1 & 0 & 0 & 0 \\
\frac{3}{11} & \frac{6}{11} & \frac{2}{11} & 0 & 0 & -1 & 0 & 0 \\
\frac{6}{11} & \frac{1}{11} & \frac{4}{11} & 0 & 0 & 0 & -1 & 0\\
\frac{7}{11} & \frac{3}{11} & \frac{1}{11}& 0 & 0 & 0 & 0 & -1
\end{array}\right)$$  

To determine the irrelevant ideal, however, we need to know the particular triangulation we are using.

Different triangulations lead to different maximal cones, and therefore a primitive collection for one fan may not be a primitive collection 
for the fan of an alternate resolution.  The differences between these irrelevant ideals means that not all crepant resolutions have the same number 
of first order deformations of the tangent sheaf.  We will focus on the existence of deformations that do not depend on the particular resolution chosen, 
and the deformations of the $G$-Hilbert scheme, a particular resolution.\\ \\

\section{Studying the Tangent Sheaf}
We now talk briefly about how to find first order deformations of the tangent sheaf of a toric variety $X$.  We recall that each first order deformation 
corresponds to an element of $\ext^1(\cT,\cT)$, where $\cT$ is the tangent bundle of $X$ \cite{Defbook}.  Let $q_i$ be the charge of $x_i$, as described 
above.  We will use as the index $i$ or $j$ for $u_1,u_2,u_3$, the vertices, and index $\alpha$ or $\beta$ for $u_4,\ldots,u_N$, the interior points.
Following \cite{Paul11} we review the argument that $\displaystyle \cD:=\bigoplus_{i=1}^3 \cO(q_i)\oplus \bigoplus_{\alpha=4}^N \cO_{D_\alpha}(q_\alpha)$ is a deformation 
of the tangent bundle.\\ \\
There is an exact sequence 
$$0\to \cO^{\oplus N-3}\to \bigoplus_{i=1}^N \cO(q_i)\to \cT\to 0$$
and we can define $\cW=\cO(q_1)\oplus\cO(q_2)\oplus \cO(q_3)$.  Since $\displaystyle\cW\hookrightarrow \bigoplus_{i=1}^N \cO(q_i)$, we 
have an induced map $f:\cW\to \cT$, where the cokernel of $f$ is the cokernel of 
$\displaystyle\cO^{\oplus N-3}\to \bigoplus_{\alpha=4}^N \cO(q_\alpha)$.  We note that the map is given by the diagonal matrix $(x_4,..,x_N)$.  
Let $D_i=V(x_i)$.  Then we know that $0\to \cO(-q_{\alpha})\overset{x_\alpha}\to \cO\to \cO_{D_\alpha}\to 0$ is  exact, 
and combining these exact sequences gives 
$$0\to \bigoplus_{i=1}^3\cO(q_i)\to \cT\to  \bigoplus_{\alpha=4}^N \cO_{D_\alpha}(q_\alpha)\to 0$$
an exact sequence.  Therefore, 
$\displaystyle \cD:=\bigoplus_{i=1}^3 \cO(q_i)\oplus \bigoplus_{\alpha=4}^N \cO_{D_\alpha}(q_\alpha)$ is a deformation 
of the tangent bundle.  So by studying $\ext^1(\cD,\cD)$, we can find all of the deformations we are looking for.  
In particular, these are only the framed deformations of $\cT$ (see \cite{Paul11}).  \\ \\

We again look at our example.  Here, we have 
$$\cD=\bigoplus_{i=1}^3 \cO(q_i)\oplus \bigoplus_{\alpha=4}^8 \cO_{D_\alpha}(q_\alpha)$$
where the $q_i$ are the columns from $\Phi$, defined in the previous section.  We note that each $\cO_{D_\alpha}$ depends on the 
resolution chosen, so $\cD$ is not independent of the triangulation.

\subsection{Deformations of the Tangent Sheaf}
We can find first order deformations of $\cD$ (and thus first order framed deformations of $\cT$) by looking at the direct sum of $\ext^1$s on each
 possible pairing of components.  There are four types of pairings that occur.  Following Aspinwall, we show they can all 
be computed using cohomology (we will do one of these computations in detail below):
\begin{enumerate}
 \item [1:]$\ext^1(\cO(q_i),\cO(q_j))=0$
 \item [2:]$\ext^1(\cO(q_i),\cO_{D_\alpha}(q_\alpha))=0$
 \item [3:]$\ext^1(\cO_{D_\alpha}(q_\alpha),\cO(q_i))=H^0(\cO_{D_\alpha}(q_i))$
 \item [4:]$\ext^1(\cO_{D_\alpha}(q_\alpha),\cO_{D_\beta}(q_\beta))=H^0(C_{\alpha\beta},N_\alpha)=H^0(\cO_{C_{\alpha\beta}}(q_\beta))$,\\
 where $N_\alpha$ is the normal bundle of $C_{\alpha\beta}$ in 
$D_\alpha$.

\end{enumerate}  
We will show (3) above, that $\ext^1(\cO_{D_\alpha}(q_\alpha),\cO(q_i))=H^0(\cO_{D_\alpha}(q_i))$.  First, recall that 
$$\ext^1(\cO_{D_\alpha}(q_\alpha),\cO(q_i))\cong \cH om^1(\cO_{D_\alpha}(q_\alpha),\cO(q_i))$$
Next, we can shift one of the exact sequences described in the previous section to find 
$$0\to \cO\overset{x_\alpha}\to \cO(q_\alpha)\to \cO_{D_\alpha}(q_\alpha)\to 0,$$  
This means that in the derived category, $\cO_{D_\alpha}(q_\alpha)$ is isomorphic to 
$\cO\overset{x_\alpha}\to \underline{\cO(q_\alpha)}$ (where the underline represents position 0)\cite{Derivedbook}.  Thus,  
$\cH om(\cO_{D_\alpha}(q_\alpha),\cO(q_i))$ 
is given by $$\underline{\cO(q_i-q_\alpha)}\overset{x_\alpha}\to \cO(q_i)$$  We recognize this as 
$\cO_{D_\alpha}(q_i)[-1]$, just a twist of the sequence described above shifted over by 1.  Therefore, we have 
$$\cH om^1(\cO_{D_\alpha}(q_\alpha),\cO(q_i))\cong \cO_{D_\alpha}(q_i)\cong H^0(\cO_{D_\alpha}(q_i))$$

%We will also show (1) above, that $\ext^1(\cO(q_i),\cO(q_j))=0$.  We have 
%$\ext^1(\cO(q_i),\cO(q_j))\cong H^1(\cO(q_j-q_i))$. We know that both $D_j$ and $D_i$  $D_j-D_i$.   
%$\cH om^1(\cO(q_i),\cO(q_j))=0$.  
The other computations can be done in a similar fashion, and can be seen in \cite{Paul11}.
We see that only two of these are non-zero, so those are the only ones we need to study.  We also see that we only need to look at global sections 
of twisted line bundles to determine what these deformations actually are.
\section{Number of Deformations of the Tangent Sheaf}
First look at the deformations that come from $\ext^1(\cO_{D_\alpha}(q_\alpha),\cO(q_i))$.  We begin by finding a basis for 
$H^0(\cO_{D_\alpha}(q_i))$.  Such a basis is given by the set of monomials in $D_\alpha$, 
whose degree in $D_\alpha$ is equal to the degree of $x_i$ in $D_\alpha$.  Since $D_\alpha=X\cap V(x_\alpha)$, 
this is the set of all Laurent monomials in $\CC(x_1,..,x_N)$ that do not contain $x_\alpha$, and which have non-negative power for all 
$x_\rho$ that share a cone with $x_\alpha$.  We know that deg($x_i$) is given by the $i$-th column of $\Phi$ (which we will 
denote $\Phi_i$).

%Recall, one column of this matrix is given by 
%$$(\nu_1^\alpha,\nu_2^\alpha,\nu_3^\alpha,0,..,-1,...,0)$$, with the $-1$ in the $\alpha$ row.
For $x_1^{c_1}\cdots x_N^{c_N}$ to have the correct multidegree, we must have $\sum c_j(\Phi_j)=\Phi_i$.  But since the 
rows of $\Phi$ and $\Phi'$ span the same complex subspace, its sufficient to find coefficients $c_j$ that make this true 
for $\Phi'$ (again, requiring $c_\alpha=0$).

We now look at this row by row. We first must have $c_1\nu_1^\alpha+c_2\nu_2^\alpha+c_3\nu_3^\alpha=\nu_i^\alpha$, where $c_j$ is 
the power of $x_j$ in the Laurent monomial.  This equation follows directly from the row corresponding to $u_\alpha$, since 
the only columns that are nonzero in this row are $x_1,x_2,x_3$, and $x_\alpha$, with $c_\alpha=0$.  
We note that the $c_j$ do not have to all be non-negative, only those where $u_j$ and $u_\alpha$ share a cone.

Once we have $(c_1,c_2,c_3)$, there is a unique $c_\beta$ for every $u_\beta$ so that 
$c_1\nu_1^\beta+c_2\nu_2^\beta+c_3\nu_3^\beta-c_\beta=\nu_i^\beta$.  Thus we have $c_\beta$ as the power of
 $x_\beta$, in order to have the correct multidegree.  Some of the $c_\beta$ can also be negative, though which
negative powers are allowed depends on the specific resolution we chose.  This leads to the following proposition.

\begin{proposition}
Every deformation from $\ext^1(\cO_{D_\alpha}(q_\alpha),\cO(q_i))$ corresponds to some triple of integers 
$(c_1,c_2,c_3)$ satisfying $c_1\nu_1^\alpha+c_2\nu_2^\alpha+c_3\nu_3^\alpha=\nu_i^\alpha$ 
with the above non-negativity conditions. 
\end{proposition}

Next, we look at $\ext^1(\cO_{D_\alpha}(q_\alpha),\cO_{D_\beta}(q_\beta))$.  To classify these deformations, we need a 
basis for $H^0(\cO_{C_{\alpha\beta}}(q_\beta))$.  This will be the set of Laurent monomials in $\CC(x_1,..,x_N)$ that do not 
contain $x_\alpha$ or $x_\beta$, and which have non-negative power for the two $x_\rho$ that share a cone with $x_\alpha$ and $x_\beta$.  We 
proceed as above, and find that the equations that must be satisfied are $c_1\nu_1^\alpha+c_2\nu_2^\alpha+c_3\nu_3^\alpha=0$ and 
$c_1\nu_1^\beta+c_2\nu_2^\beta+c_3\nu_3^\beta=-1$.  These $c_j$ will be the power of $x_j$ in the Laurent monomial.  
The only $c_j$ that must be non-negative are those that share a cone with both $u_\alpha$ and $u_\beta$.  

Given $(c_1,c_2,c_3)$, there is a unique $c_\gamma$ for each other interior points $u_\gamma$ such that 
$c_1\nu_1^\gamma+c_2\nu_2^\gamma+c_3\nu_3^\gamma-c_\gamma=0$.  As long as $c_\gamma$ is non-negative for any $u_\gamma$ sharing 
a cone with both $u_\alpha$ and $u_\beta$, $x_j^{c_j}$ gives a basis element of $H^0(\cO_{C_{\alpha\beta}}(q_\beta))$, and thus a deformation. We can 
summarize this in another proposition. 

\begin{proposition}
Every deformation from $\ext^1(\cO_{D_\alpha}(q_\alpha),\cO_{D_\beta}(q_\beta))$ corresponds to some
 triple of integers $(c_1,c_2,c_3)$ satisfying $c_1\nu_1^\alpha+c_2\nu_2^\alpha+c_3\nu_3^\alpha=0$ and 
$c_1\nu_1^\beta+c_2\nu_2^\beta+c_3\nu_3^\beta=-1$ with the above non-negativity conditions. 
\end{proposition}

The two previous propositions then give an easy way to describe any deformation of $\cD$, and therefore the tangent sheaf.

\begin{theorem}
 Every deformation of the tangent sheaf of a crepant resolution of $X$ corresponds to a unique triple of integers $(c_1,c_2,c_3)$.
\end{theorem}

In the case of $\CC^3/\ZZ_{11}$, we now look at this triple for one deformation of each type.  The triple $(2,3,0)$ gives a Laurent monomial 
$x_1^2x_2^3x_5x_6^2x_7x_8^2$, and restricting to $D_4$ (in any resolution) this gives a basis element of 
$\ext^1(\cO_{D_4}(q_4),\cO(q_3))$.  Similarly, the triple $(3,1,-2)$ gives the Laurent monomial $x_1^3x_2x_3^{-2}x_6x_7^2x_8^3$, and restricting 
to $C_{45}$ this gives a basis element of 
 $\ext^1(\cO_{D_4}(q_4),\cO_{D_5}(q_5))$ as long as $x_3=1$ on $C_{45}$.  

\subsection{Non-Isolated Singularities}
If there is a non-isolated singularity, that means that there is some lattice point $u_{\alpha}$ on the boundary of the convex hull of $(u_1,u_2,u_3)$.  
This then implies that $\nu_i^\alpha=0$ for some $i$.  But when we take this together with the propositions proved above, it means that there are 
infinitely many elements of $\ext^1(\cO_{D_\alpha}(q_\alpha),\cO(q_j))$ for $j\neq i$.  In particular, the monomials $x_i^nx_j$ will have the proper 
multidegree for all integers $n$, in any crepant resolution.  Since the main purpose of this paper is counting the number of deformations, we limit 
our focus to resolutions of orbifolds with isolated singularities.

%Now, assume that all 3 of the $c_j$ are greater than or equal to $0$.  That is, each of $x_1,x_2,x_3$ occurs with non-negative 
%multiplicity before restricting to $D_\alpha$.  This will give an element of $H^0(\cO_{D_4}(q_3))$ regardless of the resolution. 
%Either $c_3=1$, and the others are both $0$, or $c_3=0$.  
%If $c_3=1$, then all of the $c_\beta$ are $0$, and the monomial is $x_3$.\\  
%If $c_3=0$, then $c_1\nu_1^4+c_2\nu_2^4=\nu_3^4$.  
%\begin{lemma}
% If $c_1,c_2,c_3\geq 0$, then every $c_\beta\geq 0$.
%\end{lemma}

%However, this implies that 
%$c_1\nu_1^\beta+c_2\nu_2^\beta\geq\nu_3^\beta$ for all $\beta\in\{4,5,..,N\}$.  This means all $c_\beta\geq 0$, and so 
%the Laurent monomial in $\CC(x_1,..,x_N)$ is an actual monomial.  
%Therefore, any triple $(c_1,c_2,c_3)$ of non-negative 
%integers that satisfies $c_1\nu_1^4+c_2\nu_2^4+c_3\nu_3^4=\nu_3^4$ will give an element of  $H^0(\cO_{D_4}(q_3))$
%for any resolution.  Counting the number of triples of this type gives the lower bound for the 
%dimension of $\ext^1(\cO_{D_4}(q_4),\cO(q_3))$\\ \\
%Say $c_1\nu_1^\alpha+c_2\nu_2^\alpha+c_3\nu_3^\alpha=\nu_i^\alpha$.  Then we know that $c_1\nu_1^\beta

%In order to prove the second part of the theorem, we will need to use 

\section{Computing Singlets}
We will now bring in some methods from physics, which will help us to count the triples that satisfy the constraints from the previous section.  
Let $G=\ZZ_r$ acting on $\CC^3$, with G generated by $g=(\zeta_r^{a_1},\zeta_r^{a_2},\zeta_r^{a_3})$, as before.
  We associate a ``twisted sector'' with each element of $G$, where the $j/r$ twisted sector refers to $g^j$.  
We begin by studying the $j/r$-twisted sector.  We note that this means we twist by 
$$(\exp(2\pi i\nu_1),\exp(2\pi i \nu_2),\exp(2\pi i \nu_3))$$ where we define $\nu_i$ to be a rational number 
such that $r\nu_i\equiv ja_i \mod r$ and $0\leq \nu_i<1$.  We will show in the following section that these $\nu_i$ are 
in fact the same as the $\nu^\alpha_i$ found in $\Phi'$ previously, for an appropriate choice of $\alpha$.  We define
$$\widetilde{\nu}_i:=\left\{
     \begin{array}{lr}
       \nu_i-\frac{1}{2} & : \nu_i\leq\frac{1}{2}\\
			&	\\
       \nu_i-\frac{3}{2} & : \nu_i>\frac{1}{2} \\
     \end{array}
   \right.
$$
This follows from the physics, and is explained in more detail in \cite{AMP}.  We define the energy of the twisted vacuum to be 
$$E=\frac{1}{2}\sum_{i=1}^3(\nu_i(1-\nu_i)+\widetilde{\nu}_i(1+\widetilde{\nu}_i))-\frac{5}{8}$$
with charge given by $q=-\frac{3}{2}-\sum \widetilde{\nu}_i$.  For physical reasons we will not go into, we only look at
 the cases where  $\sum\nu_i=1$.  We can find all the singlets for a given twisted sector by using a generating function.  
The number of singlets is the coefficient of $q^0z^0$ in 
$$q^Ez^{-\frac{3}{2}-\sum_i\widetilde{\nu}_i}\prod_i\frac{(1+q^{1+\widetilde{\nu}_i}z^{-1})(1+q^{-\widetilde{\nu}_i}z)}{(1-q^\nu_i)(1-q^{1-\nu_i})}$$
By analyzing this partition function, we are able to arrive at a much simpler condition for counting the number of singlets 
for each twisted sector.  

\begin{proposition}
 Every singlet in a twisted sector with $(\nu_1,\nu_2,\nu_3)$ corresponds to one of the following two cases:
\begin{enumerate}

\item[Case]1: Determined by a unique non-negative triple of integers $(c_1,c_2,c_3)$ such that 
$c_1\nu_1+c_2\nu_2+c_3\nu_3=\nu_i$ for some $\nu_i$.

\item[Case]2: Determined by a unique triple of integers $(c_1,c_2,c_3)$ where exactly one $c_i<-1$, and
$c_1\nu_1+c_2\nu_2+c_3\nu_3=c_i+1$.
\end{enumerate}

\end{proposition}

\section{Twisted Sectors and Interior Points}
\subsection{Action by Generator}
We now take another look at $\CC^3/\ZZ_r$, with group action given by $(\zeta_r^{a_1},\zeta_r^{a_2},\zeta_r^{a_3})$. We 
focus on the case of isolated singularities, where we always have $a_1=1,$ $a_2=a$, and $a_3=b$.  
We recall the fan in $\RR^3$ associated to this variety has one dimensional cones with minimal generators 
$(r,-a,-b),(0,1,0),$ and $(0,0,1)$.  
%We are interested in resolving $\CC^3/G$, and getting a smooth fan.  
We know that there are $r-1$ possible twisted sectors, corresponding to 
twisting by $\frac{1}{r},\frac{2}{r},\ldots,\frac{r-1}{r}$.  Looking at the $\frac{i}{r}-$twisted sector,
we have $\nu_1=\frac{i}{r}$, 
$r\nu_2\equiv ia \mod r$, and $r\nu_3\equiv ib \mod r$, with $0\leq \nu_j<1$.  
As in the previous section, we only are interested in the case where $\sum \nu_i=1$.  For each of these twisted sectors, 
look at the point $(\nu_1r, \nu_2-a\nu_1, \nu_3-b\nu_1)$.  All points of this type are in the convex hull of the three minimal 
generators, and these are in fact all lattice points in the convex hull.  
This follows from a result of Ito and Reid, relating crepant exceptional divisors of the resolution of an orbifold $\CC^3/G$ to 
the conjugacy classes of $G$ \cite{IR}.  Using the language of our paper, this result is stated below.  
\begin{proposition}
 There is a one to one correspondence between the twisted sectors and the interior lattice points of the 
cone associated to $\CC^3/\ZZ_r$.
\end{proposition}

We note that the interior point with first coordinate $i$ corresponds to the $\frac{i}{r}-$twisted sector.  If the singularity is not isolated, 
the above proposition still holds.  In that case, we find the interior lattice points (that is, minimal generators for the irreducible exceptional divisors) 
 have coordinates given by $(\nu_1r, a_1\nu_2-a_2\nu_1, a_1\nu_3-a_3\nu_1)$.

\subsection{Relation to Charge Matrix}
We have now seen that each interior points $u_\alpha$ corresponds to a particular twisted sector, and that this twisted 
sector has a triple of rational numbers $(\nu_1,\nu_2,\nu_3)$ associated with it.  
We call these $(\nu_1^\alpha,\nu_2^\alpha,\nu_3^\alpha)$, for $u_\alpha$ and have 
$$u_\alpha=(\nu^\alpha_1r,-\nu^\alpha_1a+\nu^\alpha_2,-\nu^\alpha_1b+\nu^\alpha_3)$$
Here are those triples for our running example:
\begin{center}
\begin{tabular}{|c|c|c|c|}
\hline
$u_\alpha$		& $\nu^{\alpha}_1$	   & $\nu^{\alpha}_2$		& $\nu^{\alpha}_3$	\\ \hline
$u_4$			& $\frac{1}{11}$	   & $\frac{2}{11}$		& $\frac{8}{11}$       	\\ \hline
$u_5$			& $\frac{2}{11}$	   & $\frac{4}{11}$		& $\frac{5}{11}$       	\\ \hline
$u_6$			& $\frac{3}{11}$	   & $\frac{6}{11}$		& $\frac{2}{11}$       	\\ \hline
$u_7$			& $\frac{6}{11}$	   & $\frac{1}{11}$		& $\frac{4}{11}$       	\\ \hline
$u_8$			& $\frac{7}{11}$	   & $\frac{3}{11}$		& $\frac{1}{11}$       	\\ \hline
\end{tabular}\\[5pt]  
\end{center}
Computing the cokernel of $A$, using these coordinates for the $u_\alpha$, we can show these $\nu_i^\alpha$ are the same 
as the variables with the same name used in the matrix $\Phi'$ previously.  

\section{Counting Deformations of the Tangent Sheaf}
We saw earlier that every deformation of the tangent sheaf corresponds to a unique triple of integers $(c_1,c_2,c_3)$, and that each such triple 
satisfies equations involving the $\nu_i^\alpha$, and certain negativity conditions.  In particular, which $x_i$ or $x_\alpha$ were 
allowed to have negative power was dependent on the particular resolution chosen.  But this implies that any any actual monomial
 of $\CC(x_1,..,x_N)$ (with all powers $\geq 0$) will correspond to a deformation for any resolution.  We will now find all such monomials, by 
showing that any triple $(c_1,c_2,c_3)$ of non-negative integers, satisfying the appropriate equations, will have non-negative powers for all $x_\alpha$.

\begin{theorem}
\label{bigtheorem}
 $\ext^1(\cO_{D_\alpha}(q_\alpha),\cO(q_i))$ has a minimum number of elements for any crepant resolution of $\CC^3/G$, and the $G$-Hilbert scheme 
achieves that minimum.
\end{theorem}
We prove the first part immediately, and the second part later in the paper, when we study the $G$-Hilbert scheme.  For ease of notation, we will 
assume $i=3$ and $\alpha=4$.  Let $c_1,c_2,c_3\geq 0$, and $c_1\nu_1^4+c_2\nu_2^4+c_3\nu_3^4=\nu_3^4$.  The key step in this proof uses the following lemma.
\begin{lemma}
\label{repeatlemma}
 If $c_1,c_2,c_3\geq 0$, then every $c_\beta\geq 0$.
\end{lemma}
Combining this lemma with Proposition 4.1, any triple $(c_1,c_2,c_3)$ of non-negative 
integers that satisfies $c_1\nu_1^4+c_2\nu_2^4+c_3\nu_3^4=\nu_3^4$ will give an element of  $H^0(\cO_{D_4}(q_3))$
for any resolution.  Counting the number of triples of this type thus gives the lower bound for the 
dimension of $\ext^1(\cO_{D_4}(q_4),\cO(q_3))$.  We now prove the lemma:\\

Let the group action be given by $(a_1,a_2,a_3)$, and say $c_1\nu_1^\alpha+c_2\nu_2^\alpha+c_3\nu_3^\alpha=\nu_i^\alpha$.  
Without loss of generality, we can assume $i=3$.  We look at all triples that have each of $c_1,c_2,c_3\geq 0$.  
There is always one such triple, $(0,0,1)$.  This corresponds to the Laurent monomial $x_3$, and will obviously 
always occur (all the $c_\beta$ are 0).  All the other triples must have $c_3=0$, so our equation becomes 
$c_1\nu_1^\alpha+c_2\nu_2^\alpha=\nu_3^\alpha$. We know that 
$\nu_1^\alpha+\nu_2^\alpha+\nu_3^\alpha=1$.  Therefore, we can rewrite the above as 
$(c_1+1)\nu_1^\alpha+(c_2+1)\nu_2^\alpha=1$.  Let $c_1'=c_1+1$, $c_2'=c_2+1$.

Now, we use the fact that $r\nu_1^\alpha=ia_1 \mod r$ and $r\nu_2^\alpha=ia_2 \mod r$, so we can rewrite the equation as 
$c_1'(ia_1 \mod r) + c_2'(ia_2 \mod r)= r$.  Thus, $c_1'(ia_1-n_ir)+c_2'(ia_2-m_ir)=r$, and 
$i(c_1'a_1+c_2'a_2)=r(1+n_i+m_i)$.  We note all of the constants used are integers.  Since $i$ does not 
divide $r$, it must divide $(1+n_i+m_i)$ and we thus have $c_1'a_1+c_2'a_2=rk$ where $k$ is an integer, with $k\geq 1$.
%Therefore, we must have $c_1'a_1+c_2'a_2\geq r$.\\ \\

Next we show that for any $\beta$, $c_1\nu_1^\beta+c_2\nu_2^\beta\geq \nu_3^\beta$.  This is equivalent to showing 
$c_1'\nu_1^\beta+c_2'\nu_2^\beta\geq 1$, by the same argument as above.  
We again have $r\nu_1^\beta=ja_1 \mod r$, $r\nu_2^\beta=ja_2 \mod r$, for some $j$.  So if we 
multiply both sides by $r$, we get another equivalent condition, $c_1'(ja_1-n_jr)+c_2'(ja_2-m_jr)\geq r$.
The left hand side can then be written as $j(c_1'a_1+c_2'a_2)-r(c_1n_j+c_2m_j)$.  But this is 
just $r(jk-c_1n_j+c_2m_j)$.  Therefore, we see that the left hand side is in fact divisible by $r$.  We know that 
$c_1'(ja_1 \mod r) + c_2'(ja_2 \mod r)>0$, and therefore it must be at least $r$.  So we have shown 
$c_1'(ja_1-n_jr)+c_2'(ja_2-m_jr)\geq r$, and equivalently, $c_1\nu_1^\beta+c_2\nu_2^\beta\geq \nu_3^\beta$.  This implies that all the $c_\beta$ 
are non-negative, since $c_1\nu_1^\beta+c_2\nu_2^\beta-c_\beta=\nu_3^\beta$.  
$\square$ \\

We also note that the conditions necessary for the theorem were a triple of integers $(c_1,c_2,c_3)$ 
with $c_1,c_2,c_3\geq 0$ and $c_1\nu_1^\alpha+c_2\nu_2^\alpha+c_3\nu_3^\alpha=\nu_i^\alpha$.  Such a triple also 
uniquely determined all singlets of case 1.  This leads to the following:\\
\begin{theorem}
 Every singlet of case 1 gives a distinct first order deformation of the tangent sheaf of any crepant resolution of $\CC^3/G$.
\end{theorem}

\subsection{Other Singlets}

Unfortunately, Theorem 7.1 does not extend to singlets of case 2.  Recall, such a singlet is identified by a triple of integers 
$(c_1,c_2,c_3)$ where exactly one $c_i<-1$, and $c_1\nu^\alpha_1+c_2\nu^\alpha_2+c_3\nu^\alpha_3=c_i+1$. We find that whether such a singlet corresponds 
to a deformation of $\cT$ depends on the particular triangulation chosen.  However, there is still a relationship to deformations.

\begin{proposition}
 If $(c_1,c_2,c_3)$ identify a singlet, then they also uniquely identify a deformation of the tangent sheaf of some particular resolution of $\CC^3/G$.
\end{proposition}
We have already shown this proposition for singlets of case 1.  We now show this for singlets of case 2.
Say we have a triple of integers as described above.  For ease of reading, we assume that $i=3$, so 
$c_1\nu^\alpha_1+c_2\nu^\alpha_2+c_3\nu^\alpha_3=c_3+1$, where $c_3$ is negative.  Since we know that $\nu^\alpha_3=1-\nu^\alpha_1-\nu^\alpha_2$, 
we can rewrite this as $(c_1-c_3)\nu_1^\alpha+ (c_2-c_3)\nu_2^\alpha=1$.  We also know the coordinates of $u_\alpha$ are 
$(\nu^\alpha_1r, a_1\nu^\alpha_2-a_2\nu^\alpha_1, a_1\nu^\alpha_3-a_3\nu^\alpha_1)$.
%$(ia_1-n_ir,n_ia_2r-m_ia_1r,1-ia_1+n_ir(1-a_2)+m_ira_1)$.
Now look at $$u_\beta=(-c_3-1)u_\alpha\qquad \text{ and } \qquad u_\gamma=(-c_3)u_\alpha.$$
 % This implies $\nu^\beta_1=(-c_3-1)\nu^\alpha_1$ and
 %$\nu^\beta_2=(-c_3-1)\nu^\alpha_2$, and similarly, $\nu^\gamma_1=-c_3\nu^\alpha_1$ and
 %$\nu^\gamma_2=-c_3\nu^\alpha_2$. 
After calculating, we find $c_1\nu_1^\gamma+c_2\nu_2^\gamma+c_3\nu_3^\gamma=0$, and 
$c_1\nu_1^\beta+c_2\nu_2^\beta+c_3\nu_3^\beta=-1$.
%$$c_1\nu_1^\gamma+c_2\nu_2^\gamma+c_3\nu_3^\gamma=(c_1-c_3)\nu_1^\gamma+(c_2-c_3)\nu_2^\gamma+c_3=$$
%$$(c_1-c_3)(-c_3)\nu_1^\alpha+(c_2-c_3)(-c_3)\nu_2^\alpha+c_3=-c_3((c_1-c_3)\nu_1^\alpha+ (c_2-c_3)\nu_2^\alpha-1)=0$$
%and 
%$$c_1\nu_1^\beta+c_2\nu_2^\beta+c_3\nu_3^\beta=(c_1-c_3)\nu_1^\beta+(c_2-c_3)\nu_2^\beta+c_3=$$
%$$(c_1-c_3)(-c_3-1)\nu_1^\alpha+(c_2-c_3)(-c_3-1)\nu_2^\alpha+c_3+1-1=(-c_3-1)((c_1-c_3)\nu_1^\alpha+ (c_2-c_3)\nu_2^\alpha-1)-1=-1$$

We note that these are exactly the equations we needed for $(c_1,c_2,c_3)$ to give a deformation from the group 
$\ext^1(\cO_{D_\alpha}(q_\alpha),\cO_{D_\beta}(q_\beta))$.$\square$

While these particular deformations do not show up for every crepant resolution of $\CC^3/G$, we note that there are always at least this many 
deformations.  This can be seen from the geometry in the following section.

This proposition leads to the question of whether there is a resolution with a one to one correspondence between first order deformations and singlets.  
We will show the answer is yes, and that the $G$-Hilbert scheme is such a resolution.  
In order to do this, we will first construct a quiver associated with the original orbifold $\CC^3/G$.

\section{Quiver Construction of a Resolution}
We will now use the two types of singlets to construct a quiver associated to the orbifold $\CC^3/G$.  This quiver will then lead to a distinguished 
triangulation of $\{u_1,..,u_N\}$, and thus a particular resolution of the orbifold.

  As a reminder, for each sector (and thus each $u_\alpha)$, we had two possible cases of singlets:
\begin{enumerate}

\item[Case]1: Determined by a unique non-negative triple of integers $(c_1,c_2,c_3)$ such that 
$c_1\nu^\alpha_1+c_2\nu^\alpha_2+c_3\nu^\alpha_3=\nu^\alpha_i$ for some $\nu^\alpha_i$.

\item[Case]2: Determined by a unique triple of integers $(c_1,c_2,c_3)$ where exactly one $c_i<-1$, and
$c_1\nu^\alpha_1+c_2\nu^\alpha_2+c_3\nu^\alpha_3=c_i+1$.
\end{enumerate}

For each singlet of case 1, include an arrow from $u_\alpha$ to $u_i$.  We note that there is always at least one arrow to each vertex, 
since $(1,0,0),(0,1,0),$ and $(0,0,1)$ all meet the criteria.  Below we have the quiver at this stage, for the $\CC^3/\ZZ_{11}$ example we have been 
working with.
\begin{center}
\includegraphics[height=2in]{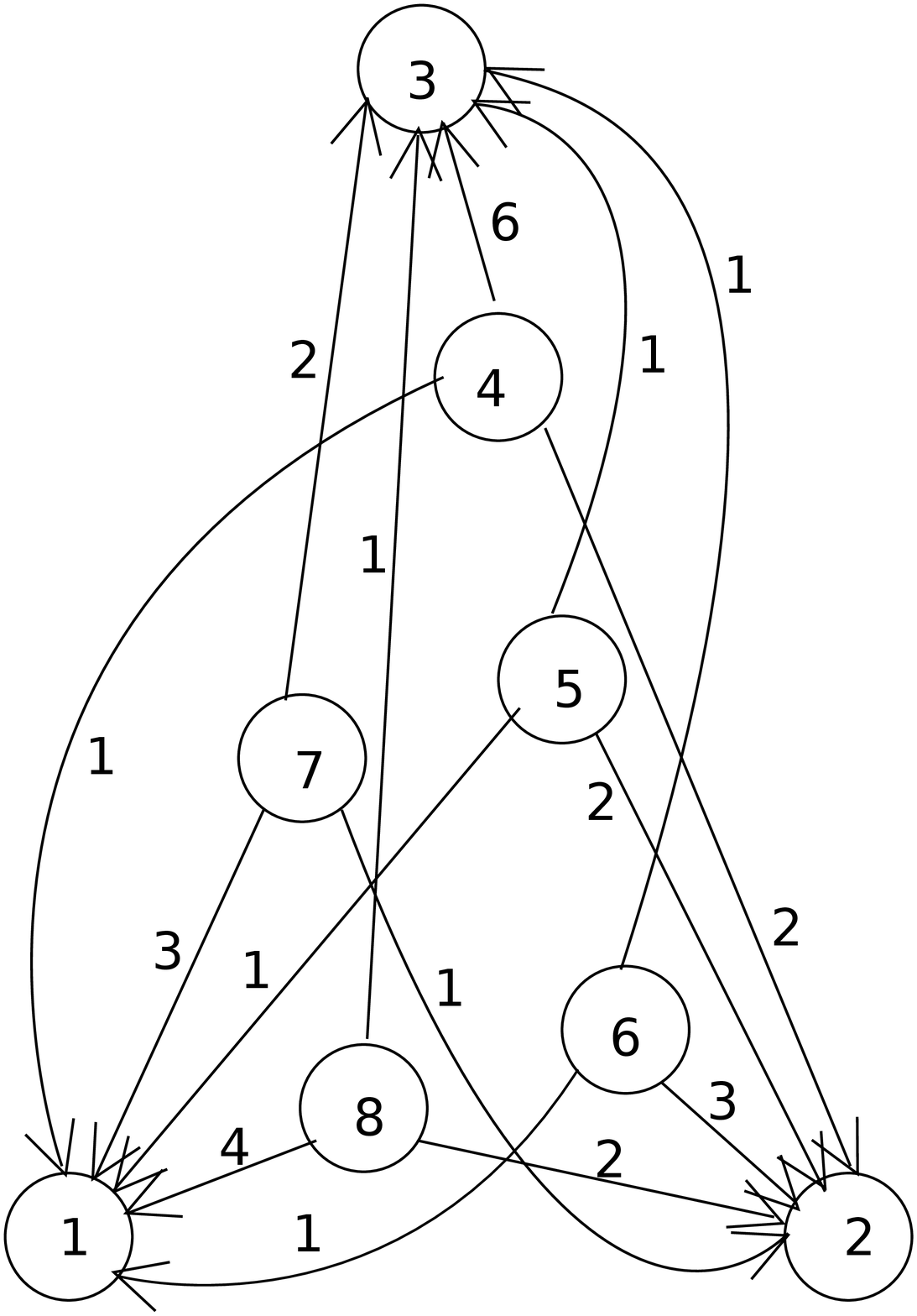}\\
Figure 8.1:  First step of quiver for $\CC^3/\ZZ_{11}$\\
\end{center}

Now, for each singlet of case 2 we know exactly one $c_i$ is negative.  Extend the line from 
$u_i$ through $u_\alpha$ until it passes through $-c_i$ nodes total, in addition to $u_i$.  For each singlet with this value for $c_i$, include an 
arrow from this node to the node one closer to $u_i$.  We note that if there is a singlet with $c_i=-n$, then there will be singlets with 
$c_i=-m$ for all $m<n$, so there are arrows between all nodes along this line, pointing towards $u_i$.  
Let's look at our example again, concentrating first on $u_4$.  In this case, we had $\nu_1^4=\frac{1}{11}$, $\nu_2^4=\frac{2}{11}$, and 
$\nu_3^4=\frac{8}{11}$.  The only triples $(c_1,c_2,c_3)$ that satisfy the condition for case 2 are 
$(5,0,-2),(3,1,-2),(1,2,-2),(1,0,-3),$ and $(0,1,-3)$.  This means we must extend the line from $u_3$ to $u_4$ through two more nodes 
($u_5$ and $u_6$), with arrows of multiplicity 3 from $u_5$ to $u_4$, and multiplicity 2 from $u_6$ to $u_5$.  These arrows, along with 
those corresponding to the other $u_\alpha$, are included in the final quiver.
   \begin{center}
\includegraphics[height=2in]{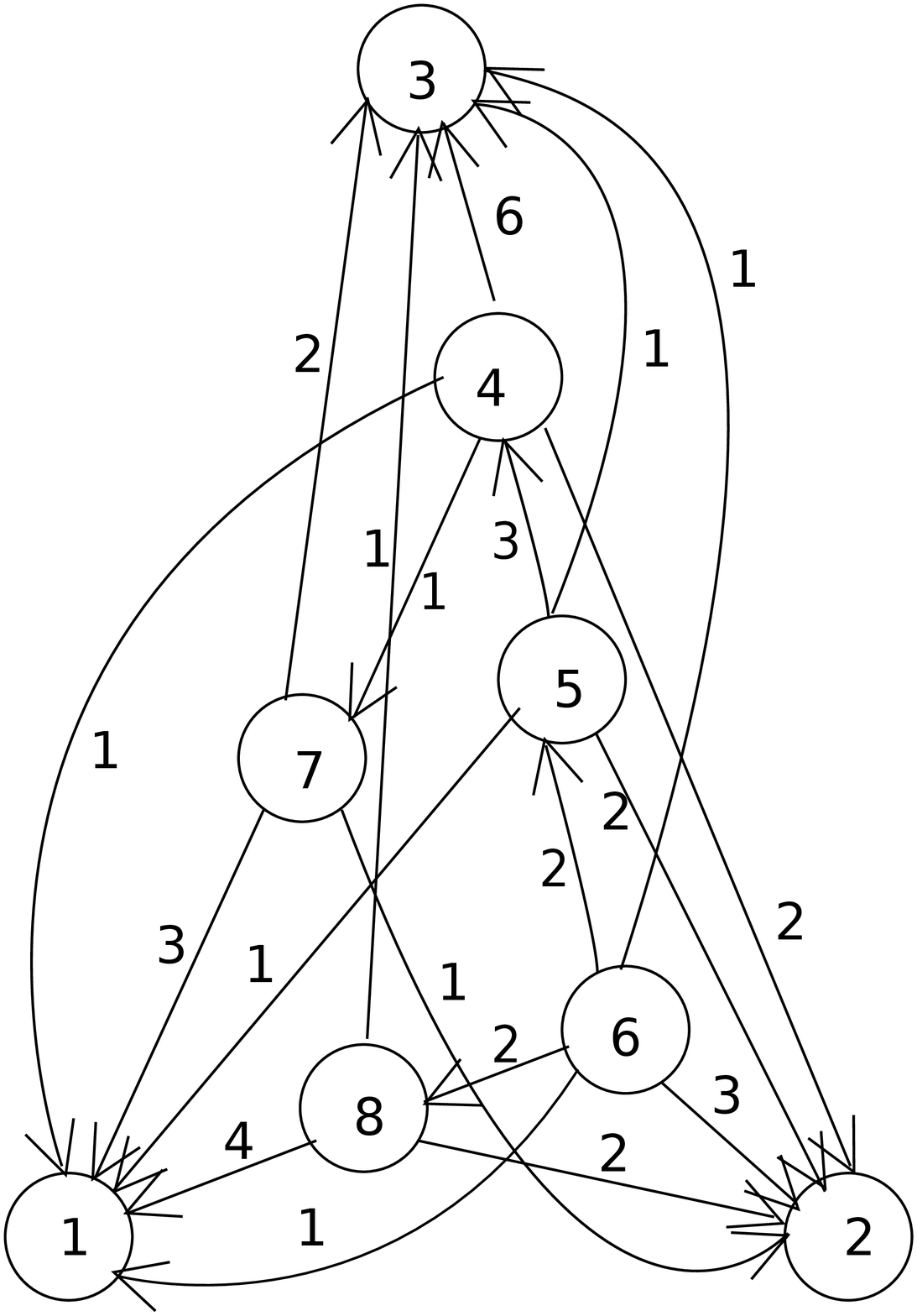}\\
Figure 8.2:  $\CC^3/\ZZ_{11}$ quiver\\
\end{center}

\subsection{Choosing a Triangulation}
We now use this quiver to construct a triangulation of the point set $\{u_1,..,u_N\}$.  
First, we suppress all arrows to vertices of multiplicity one.  Since they always exist we don't lose any information by not 
including them, but now the quiver is planar.  Next, replace every arrow with an edge and add in the 3 boundary edges, connecting the vertices to each other.  Finally, 
take a regular tessellation of any remaining triangles that are not smooth.  In this way, we get a smooth triangulation of the triangle with 
vertices $\{u_1,u_2,u_3\}$.  

Here is the triangulation for $\CC^3/\ZZ_{11}$.  We include the multiplicities 
the arrows previously had, since they will prove useful later on.

  \begin{center}
\includegraphics[height=2in]{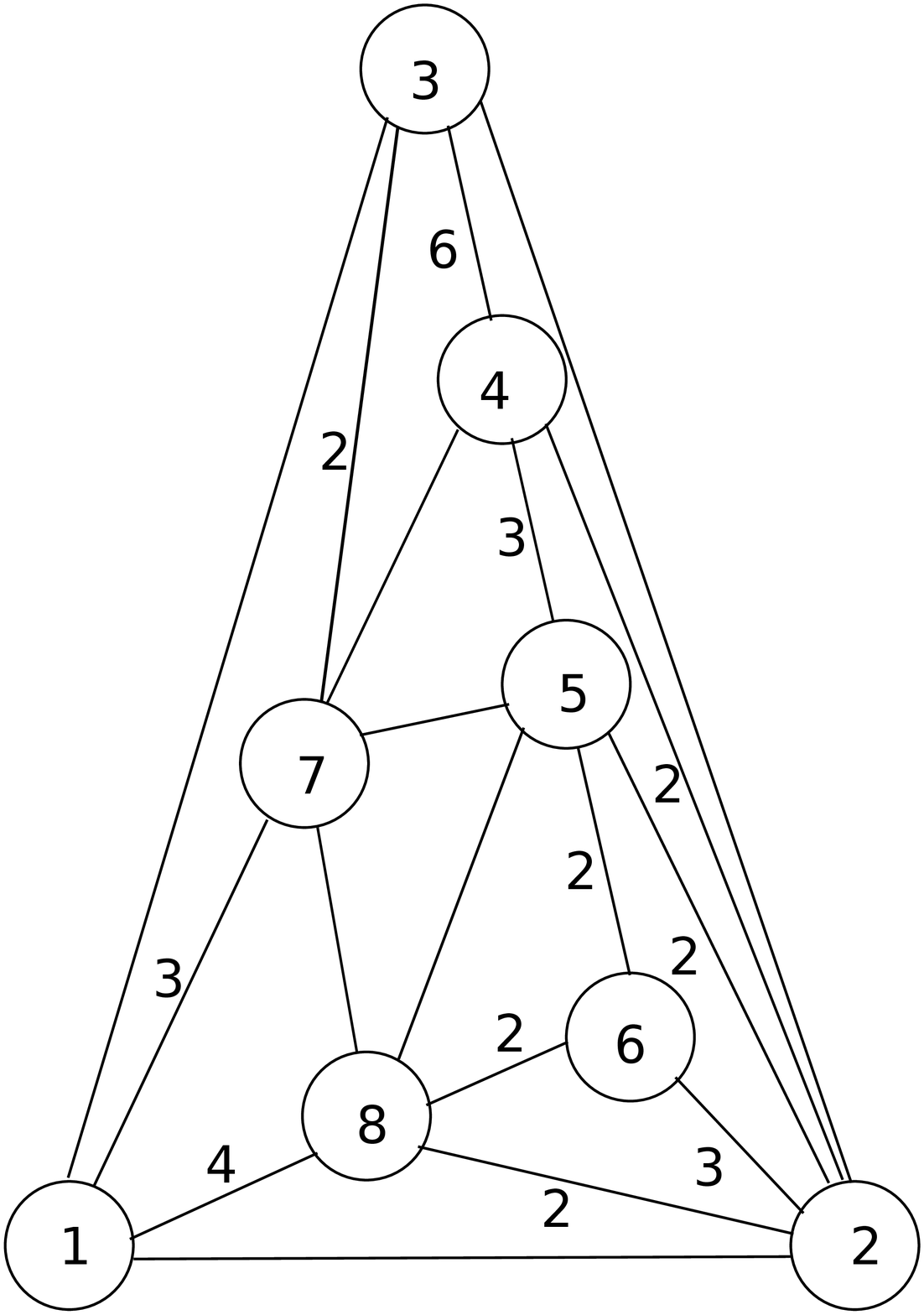}\\
Figure 8.3:  $\CC^3/\ZZ_{11}$ Triangulation.\\
\end{center}

If we now take a fan over this triangulation, we get a crepant resolution of the orbifold $\CC^3/G$.  
We will show that this resolution is actually the $G$-Hilbert scheme.  To do this, we will expand on the knockout method 
for constructing the fan of the G-Hilbert scheme described in \cite{CR}, with a slight modification.

\section{$G$-Hilbert Scheme}
\subsection{Jung-Hirzebruch Continued Fraction}
The $G$-Hilbert scheme is a particular resolution of $\CC^3/G$, which parameterizes the $0$-dimensional subschemes of $\CC^n$ that are invariant 
under the group action, have length $|G|$, and have global sections that give a regular representation of $G$ \cite{IN}.  Since this scheme is a toric 
variety, we know that there is a fan associated with it.  \cite{CR} gives us a combinatorial method for constructing this associated fan, which we now review.

The first step in constructing the fan of the $G$-Hilbert scheme is finding the Jung-Hirzebruch continued fraction associated to each vertex.  We 
include here a quick review of how the continued fraction is computed.  
For each vertex ($u_1,u_2$, and $u_3$), we look at the two dimensional cone whose origin is the vertex and whose walls are the rays 
pointing to the other two vertices.  This cone will likely be singular, but as a toric surface there's a well known method for resolving it using the 
continued fraction.  Choose a basis so that one ray has minimal generator $(1,0)$ and the other has generator $(-c,r)$.  We can then write 
$r/c$ as a continued fraction, $r/c=[[b_1,..,b_m]]$.  This means that $$\frac{r}{c}=b_1-\frac{1}{b_2-\frac{1}{b_3-...}}$$  We can find the 
$b_i$ through the equations 
$r=c(b_1)-d_1, c=d_1(b_2)-d_2, d_1=d_2(b_3)-d_3,\ldots, d_{m-2}=d_{m-1}b_m-d_m$.  There are two more values associated with this continued fraction, 
which we will define recursively.  We define $P_i$ by setting $P_0=1,P_1=b_1$, and $P_i=b_iP_{i-1}-P_{i-2}$.  Similarly, we define $Q_i$ by setting 
$Q_0=0,Q_1=1$, and $Q_i=b_iQ_{i-1}-Q_{i-2}$.  In particular, note that this implies $P_i=b_i(P_{i-1}+P_{i+1})$ and $Q_i=b_i(Q_{i-1}+Q_{i+1})$
  
For each vertex, $r$ is the order of the underlying group $G$, 
and the particular $c$ can be found through the group action.  We then have $m$ rays from the vertex to interior points, and label each ray with the number 
$b_i$, which we call its strength.  The coordinates of the interior points, in the basis we chose earlier, will be $(P_i,Q_i)$.   
We do this for the other two vertices as well, so each vertex has a number of labeled rays to particular interior points.

Let's look at our example, $\CC^3/\ZZ_{11}$.  At vertex $u_3$, we have $11/2=[[6,2]]$ since $11=2(6)-1$, and $2=1(2)-0$.  So the 
edges have strength $6,2$ respectively.  At $u_2$ the continued fraction is $11/7=[[2,3,2,2]]$, and at $u_1$ it is $11/4=[[3,4]]$.

\begin{center}
\includegraphics[height=2in]{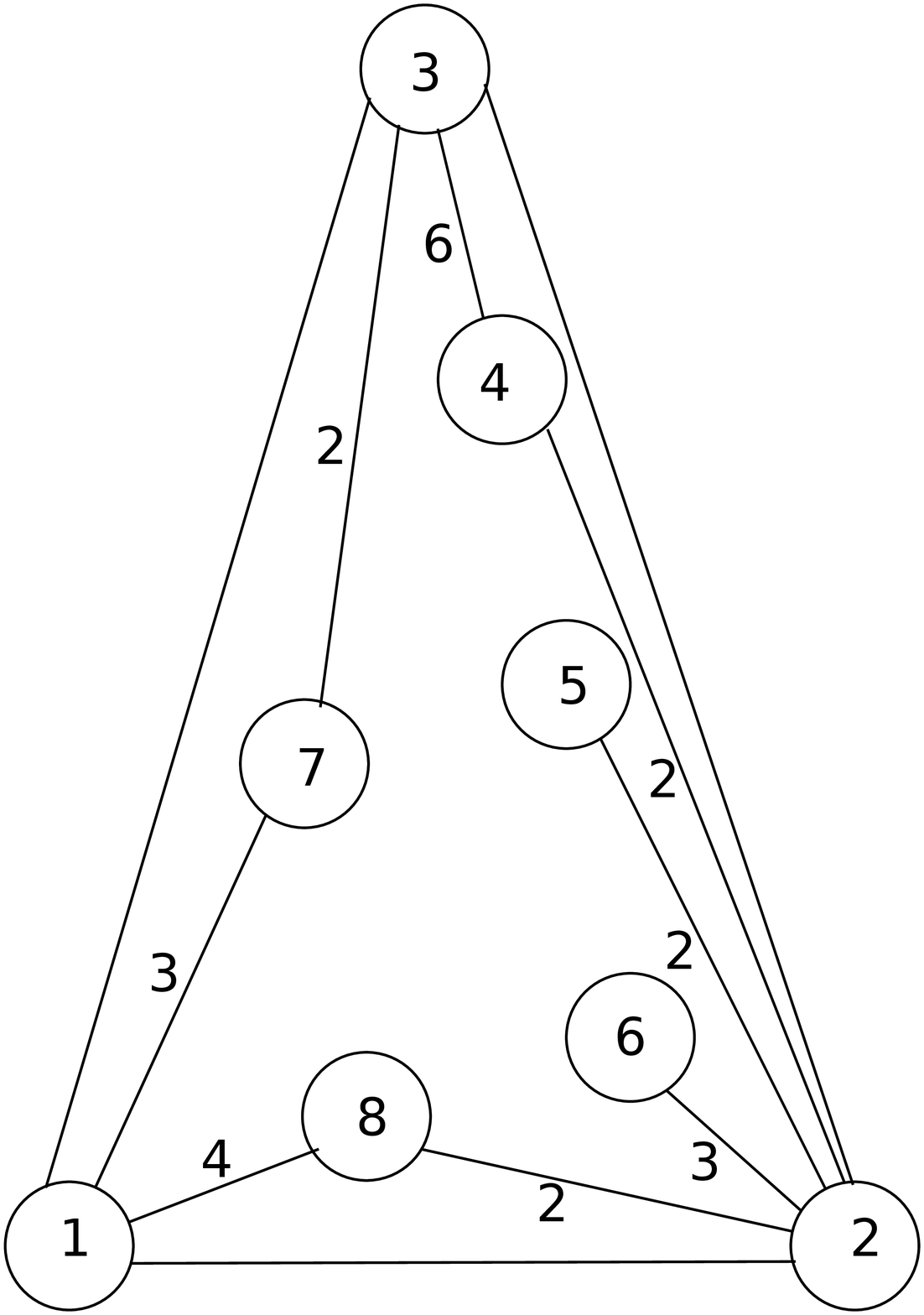}\\
Figure 9.1:  First Step in Constructing $G$-Hilbert Scheme\\
\end{center}

\subsection{Smooth Triangulation}
  Once we have the lines from vertices, extend each line further inside of the triangle.  As it goes through interior points, reduce the strength by 1 
for every other line that goes through the same interior point.  If a line would ever have 1 or 0 strength, do not include it.  Finally, take a 
regular tessellation of any remaining regular triangles.  The fan over this triangulation gives the $G$-Hilbert scheme.  
(Thus far, this is exactly the knockout method from \cite{CR}).

We now modify this method slightly by subtracting an additional 1 the first time an edge goes between 2 interior points.  This results in 
all the numbers on interior edges being 1 less than they are determined to be in \cite{CR}, and means it is now possible to have interior edges of 
strength 1.

With this modification, here is the triangulation (with strengths) for the $G$-Hilbert scheme of $\CC^3/\ZZ_{11}$.
 \begin{center}
\includegraphics[height=2in]{C3Z11picturecok.eps}\\
Figure 9.2:  Triangulation for $G$-Hilbert Scheme of $\CC^3/\ZZ_{11}$.\\
\end{center}
We see that this is the same triangulation as was found in the previous section, constructed using the singlets.  Not only is the triangulation the same, 
but the numbers along each edge are also the same.  The remainder of the paper will focus on proving that this is always the case, and understanding 
what this tells us about the $(0,2)-$McKay correspondence.

\section{Deformations of the Tangent Sheaf of the G-Hilbert Scheme}
Every framed deformation of the tangent sheaf $\cT$ corresponds to an element of 
$\ext^1(\cO_{D_\alpha}(q_\alpha),\cO(q_i))$ or an element of 
$\ext^1(\cO_{D_\alpha}(q_\alpha),\cO_{D_\beta}(q_\beta))$.  The above method of constructing 
the Hilbert scheme then gives the dimension for each of these groups, as follows.  The strength of the edge between $u_\alpha$ and $u_i$ gives the 
dimension of $\ext^1(\cO_{D_\alpha}(q_\alpha),\cO(q_i))$. The strength of the edge between $u_\alpha$ and $u_\beta$ gives the 
dimension of $\ext^1(\cO_{D_\alpha}(q_\alpha),\cO_{D_\beta}(q_\beta))$, where $u_\beta$ is closer to the vertex that gave this edge
 (see previous section).  We will now summarize and prove this relationship.
\begin{theorem}
 The numbers on edges as described above will be exactly the dimension of the $\ext^1$ group corresponding to the given vertices, and therefore count 
the number of framed deformations.
\end{theorem}

\subsection{Proof of Theorem}
We first show that the numbers on edges that go from $u_\alpha$ to $u_i$ match up with the dimension of $\ext^1(\cO_{D_\alpha}(q_\alpha),\cO(q_i))$. 

Call the $u_i$ vertex $(0,0)$, where corresponding coordinates for $u_\alpha$ are $(x,y)$, and the two neighboring vertices are $(a,b)$ and $(c,d)$.  
By definition of continued fraction, we know that $a+c=Nx$ and $b+d=Ny$, where $N$ is the strength of the edge.  Now let's move the 
origin to $(x,y)$, as in \cite{Paul11}.  In this case, we know the dimension is given by $x_2y_3-x_3y_2+2$.  Now the coordinates are $(-x,-y)$ for top one, $(a-x,b-y)$ and $(c-x,d-y)$.  We compute $x_2y_3-x_3y_2$, 
and find $(a-x)(d-y)-(b-y)(c-x)$.  We know $x=\frac{a+c}{N}$ and $y=\frac{b+d}{N}$, so this becomes 
$$\frac{1}{N^2}(((N-1)a-c)((N-1)d-b)-((N-1)c-a)((N-1)b-d))$$  which is just $\frac{1}{N}(N-2)(ad-bc)$.  But we note that 
$(Nx,Ny)\times (c,d)=((a,b)+(c,d))\times (c,d)=ad-bc$, and $(Nx,Ny)\times (c,d)=N(x,y)\times (c,d)=N$, so  $x_2y_3-x_3y_2=N-2$.  Therefore, the 
dimension of $\ext^1(\cO_{D_\alpha}(q_\alpha),\cO(q_i))$ is $N-2+2=N$.
% since this triangulation is smooth and these rays are next to each other.  So in notation of paper, $m=N-2$.
Thus, we have shown that the strength $N$, is exactly the dimension of the Ext group for edges to vertices.\\

Now look at edges that do not go to vertices.  Begin by setting one of the nodes on the edge, $u_\alpha$, as the origin.  We then have 
coordinates for the other three nodes in this system, $u_i=(a,b)$, $u_j=(c,d)$ and $u_\beta=(x,y)$.  Since we are using a smooth triangulation, 
we know that $(a,b)+(c,d)=N(x,y)$ for some $N$.  If this number is 0 or negative, then instead set $u_\beta$ as the origin, which will give an 
integer $N\geq 1$.  Note that the positive $N$ satisfying this equation came from using the Hirzebruch-Jung resolution and contracting, and 
therefore this positive $N$ is the (unmodified) strength as originally found in \cite{CR}.  
This also means that in the modified version, the strength of the edge is $N-1$.
%Thus, $(Nx,Ny)\times (c,d)=((a,b)+(c,d))\times (c,d)=(a,b)\times (c,d)=ad-bc$.  
%But we also have $(Nx,Ny)\times (c,d)=N(x,y)\times (c,d)=N$, since this triangulation is smooth and these rays are next to each other.  We therefore 
%we have $N=ad-bc$.  If we had instead set $u_\beta$ as the origin, we get 
%$((a,b)-(x,y))+((c,d)-(x,y))=N'(-x,-y)$.  Using  $(a,b)+(c,d)=N(x,y)$, we get $(N-2)(x,y)=-N'(x,y)$, and therefore $N'=2-N$.  We just refer to 
%the larger as $N$.  
As in \cite{Paul11}, the dimension of $\ext^1(\cO_{D_\alpha}(q_\alpha),\cO_{D_\beta}(q_\beta))$ is given by  $x_2y_3-x_3y+1$.  After computing as in 
the previous case, we find this is equal to $N-1$.  Therefore, the numbers on each interior edge as described in the previous method do give 
the dimension of $\ext^1(\cO_{D_\alpha}(q_\alpha),\cO_{D_\beta}(q_\beta)).\square$

By studying the geometry of this triangulation, we see that any other triangulation will have at least the same total number of interior deformations.  
This is because the only edges we can change are those with strength 0, and therefore are part of a regular tessellation.  
Any flop will introduce new deformations on bordering strength 0 edges, and as there are always at least two of those, the total number of interior 
deformations cannot decrease.

\section{Relationship between Singlets and the \\$G$-Hilbert Scheme}

We will now prove that the triangulation given by the singlets, from section 8, is in fact the triangulation which gives the G-Hilbert scheme.  
To show the triangulation is the same, it will be enough to show that the numbers on each edge from the modified knockout method (section 9), 
and the quiver method (section 8) agree.  This is because in both cases, the triangulation is uniquely determined by these numbers.  To do so, we 
first find a relationship between the $\nu_i^\alpha$ and the $P_i$ that arise from the continued fraction.  For ease of notation, we focus on the 
continued fraction based at vertex $u_3$, and the edges to $u_3$.  In particular, this means whenever we write $P_i$, we will be referring to that 
continued fraction.  All of the results generalize to the other vertices, but writing the computations in complete 
generality is a bit unwieldy.
\subsection{Edges to a Vertex}
\begin{proposition}
 For $u_\alpha$ that appears as the $j$-th term in the Hirzebruch-Jung resolution based at $u_3$, we have 
$\nu^\alpha_1=\frac{P_j}{r}$, $\nu^\alpha_2=\frac{d_j}{r}$, $\nu^\alpha_3=\frac{r-P_j-d_j}{r}$, and $r=P_jd_{j-1}-P_{j-1}d_j$
\end{proposition}
We will now prove the above relations.  $\nu^\alpha_1=\frac{P_i}{r}$, since $P_i=(u_\alpha)_1=\nu^\alpha_1r$.  Similarly, 
we have $\nu_2=\frac{aP_i \mod r}{r}$.  We can then show by induction that $P_ia\equiv d_i \mod r$, for all $i>1$.  
%$P_1a=b_1a=r+d_1\equiv d_1 \mod r$.  We now look at $P_{k+1}$.  $P_{k+1}a=(b_{k+1}P_k-P_{k-1})a \equiv b_{k+1}d_k-d_{k-1}=d_{k+1}$.  
Since every $d_i<r$, and $0\leq \nu^\alpha_2<1$, we must actually have $\nu^\alpha_2=\frac{d_i}{r}$.  The condition on $\nu^\alpha_3$ follows from
 $\sum_i \nu^\alpha_i=1$.  
Finally, we know that $P_1d_0-P_0d_1=b_1a-d_1=r$, by the first relation of the Hirzebruch-Jung resolution.  
We note we also have $d_{n-1}=b_{n+1}d_n-d_{n+1}$.  
Therefore $P_nd_{n-1}-P_{n-1}d_n=b_{n+1}P_nd_{n}-P_nd_{n+1}-P_{n-1}d_n=(b_{n+1}P_n-P_{n-1})d_n-P_nd_{n+1}$.
But $b_{n+1}P_n-P_{n-1}=P_{n+1}$ by the recursion relation, so this is $P_{n+1}d_n-P_nd_{n+1}$.  Since we have shown all of these combinations are the 
same, they are thus all equal to $r$.  We see that we can then rewrite $\nu_3$ as $\frac{P_i(d_{i-1}-1)-(P_{i-1}+1)d_i}{r}.\square$\\

%Note: There is a similar relationship involving the $P_j$ arising from the continued fraction based at the other vertices.  In those cases, the 
%formulae for the $\nu^\alpha_i$s are permuted.  For example, when the continued fraction is based at $u_1$, we have 
%$$\nu^\alpha_1=\frac{r-P_j-d_j}{r},\qquad \nu^\alpha_2=\frac{P_j}{r},\qquad \nu^\alpha_3=\frac{d_j}{r}$$

Now, recall the way we determine the number of arrows from $u_\alpha$ to $u_3$ in the singlet construction is by finding all non-negative triples with 
$c_1\nu_1+c_2\nu_2+c_3\nu_3=\nu_3$.  We always have $c_1=c_2=0, c_3=1$.  We want to count how many triples there are so that $c_1\nu_1+c_2\nu_2=\nu_3$.  
To do so, we will use the relationships we just found.  This gives $c_1\frac{P_i}{r}+c_2\frac{d_i}{r}=\frac{P_i(d_{i-1}-1)-(P_{i-1}+1)d_i}{r}$ and thus 
$$c_1P_i+c_2d_i=P_i(d_{i-1}-1)-(P_{i-1}+1)d_i$$ Now say that $c_2=P_in-(P_{i-1}+1)$, where $n$ is some integer.  We then have 
$$c_1P_i+(P_in-(P_{i-1}+1))d_i=P_i(d_{i-1}-1)-(P_{i-1}+1)d_i$$

A bit more algebra gives $c_1P_i+P_ind_i=P_i(d_{i-1}-1)$, and we can divide out by $P_i$ to get  $$c_1+nd_i=d_{i-1}-1$$
We now see that when $n=1$, this gives $c_1=d_{i-1}-1-d_i\geq 0$, since $d_i$ is decreasing sequence, and $c_2=P_i-P_{i-1}-1\geq 0$ since $P_i$ is 
increasing sequence.  Therefore, $n=1$ gives a triple that works.

Next look at $n=b_i-1$.  This gives $c_1+(b_i-1)d_i=d_{i-1}-1$.  We recall that $b_id_i=d_{i-1}+d_{i+1}$.  So this becomes 
$c_1+d_{i-1}-d_i+d_{i+1}=d_{i-1}-1$, and gives $c_1=d_i-d_{i+1}-1$.  $c_1\geq 0$, since $d_i$ is decreasing,
and $c_2=P_i(b_i-1)-(P_{i-1}+1)\geq 0$, since $b_i>1$.  So $n=b_i-1$ also gives a triple that works.

We see that $n$ can be any integer from 1 to $b_i-1$, giving $b_i-1$ triples that work.  Adding in the earlier case where $c_3=1$, we get a total of 
$b_i$ singlets of this type.  But $b_i$ is exactly the strength that we determined by the Hirzebruch-Jung resolution.  Thus, 
for every edge from an interior point to the vertex $x_3$, the two methods agree.

Incidentally, this also proves the second part of theorem $7.1$.
\begin{reptheorem}{bigtheorem}
 $\ext^1(\cO_{D_\alpha}(q_\alpha),\cO(q_i))$ has a minimum number of elements, and the $G$-Hilbert scheme achieves that minimum.
\end{reptheorem}

This is because we have just shown that for the $G$-Hilbert scheme, we have
$$|\ext^1(\cO_{D_\alpha}(q_\alpha),\cO(q_i))|=|\{(c_1,c_2,c_3)\in \ZZ^3| c_1,c_2,c_3\geq 0, c_1\nu^\alpha_1+c_2\nu^\alpha_2+c_3\nu^\alpha_3=\nu^\alpha_i|$$

\subsection{Interior Edges}
Next, look at the number on an interior edge.  By the modified rule above, we find this number by looking at the strength on the edge one closer 
to a vertex, subtracting one for each other edge that goes to the shared node, and subtracting an additional 1 when we go from 
an edge adjacent to a vertex to an interior one.  If the number ever drops to 0, we do not include that edge.\\

Recall that the way we determine what singlets occur between two interior nodes was by finding all possible $c_1,c_2\geq 0$ such that 
$$c_1\nu^\alpha_1+c_2\nu^\alpha_2=(1-c_3)\nu^\alpha_3+c_3$$ where $u_\alpha$ lies on the edge extended into a line, bordering $u_3$.  
First, let $n=-c_3$ and assume that the number of possibilities such that 
$c_1\nu^\alpha_1+c_2\nu^\alpha_2=(n+1)\nu^\alpha_3-n$ is $m$.  We want the 
number of possible $d_i$ so that $$d_1\nu^\alpha_1+d_2\nu^\alpha_2=(n+2)\nu^\alpha_3-(n+1)=(n+1)\nu^\alpha_3-n-(1-\nu^\alpha_3)$$  
By the method described in section 8, each pair $(d_1,d_2)$ gives an arrow on the ``next'' interior edge.  But since 
$1=\nu^\alpha_1+\nu^\alpha_2+\nu^\alpha_3$, we see this is equal to $(n+1)\nu^\alpha_3-n-\nu^\alpha_1-\nu^\alpha_2$.  So we have 
$(d_1+1)\nu^\alpha_1+(d_2+1)\nu^\alpha_2=(c_1)\nu^\alpha_1+c_2\nu^\alpha_2$, and see that the pairs 
$(d_1,d_2)$ that work are exactly $(c_1-1,c_2-1)$.   Note, we still have the requirement that $d_1,d_2\geq 0$.  So the pairs that are valid are exactly 
those with both $c_1,c_2\neq 0$.  If none of the triples had $c_1,c_2= 0$, the number of singlets is still $m$.  If one could be $0$, the number is $m-1$, and if both 
could be 0, the number is $m-2$.
Now, if there is a combination with $c_1=0$, that means the node lies on the extension of %which is $n$ away from $x_3$ 
an edge that comes from $u_1$.  
%We have c_2\nu_2=(n+1)\nu_3-n=(n+1)(1-\nu_1-\nu_2)-n=1-(n+1)\nu_1-(n+1)\nu_2.  So (c_2+(n+1))\nu_2=1-(n+1)\nu_1, which means we can write 
%\nu_1, or interior, using c_2s, exactly means we would have included the line.....  
Similarly, if there is a singlet with $c_2=0$, that means the node lies on the extension of an edge from $u_2$.  We therefore see that the strength of 
each edge changes exactly the same way as it does in the knockout method.  Finally, we note that when we are going from an ``exterior ''edge (to a vertex)
 to a purely interior one, the strength decreases by an additional one.  This is because the singlet corresponding to $c_1=c_2=0,c_3=1$ 
occurs in addition to those of the form $c_1\nu_1+c_2\nu_2=\nu_3$, and has no analog among the singlets of type 2.
So the number on the interior edges decreases exactly according to the modified knockout rule described in section 9.

We have therefore proven the following theorem:
\begin{theorem}
 The singlets associated to the orbifold $\CC^3/G$ can be used to construct the fan for the $G$-Hilbert scheme, along with the $\ext$-quiver, in a 
new way.
\end{theorem}

\subsection{Singlets and First Order Deformations}
We have proven that the distinguished triangulation given by the singlets is in fact the G-Hilbert scheme. In addition, by the theorems of section 7, we have also shown
 that every first order deformation of the $G$-Hilbert scheme corresponds to a particular singlet. This has a number of useful implications.  
First of all, it implies that the quiver constructed in section 6 is actually the $\ext-$quiver for the $G$-Hilbert scheme.  
That is, if we associate each summand $V_i$ of $\cD=\bigoplus_{i=1}^3 \cO(q_i)\oplus \bigoplus_{\alpha=4}^8 \cO_{D_\alpha}(q_\alpha)$ with the 
node $u_i$, there are (dim $\ext^1(V_i,V_j)$) arrows from node $i$ to node $j$.

Second, this shows the $(0,2)$-McKay correspondence always works for the $G$-Hilbert scheme.  The total number of first order deformations 
of the tangent sheaf is the same as the number found from the conformal field theory.  In addition, we have the further result that the number of 
singlets of case 1 gives the number of deformations that come from $\ext^1(\cO_{D_\alpha}(q_\alpha),\cO(q_i))$,
 while the number of singlets of case 2 gives the number of deformations coming from $\ext^1(\cO_{D_\alpha}(q_\alpha),\cO_{D_\beta}(q_\beta))$.

Also, we now have proven all parts of theorem 1.1:
\begin{reptheorem}{maintheorem}
 Among all crepant resolutions of $\CC^3/\ZZ_r$ with $r$ prime, the $G$-Hilbert scheme has minimal dimension for $\ext^1(\cT,\cT)$.  This dimension exactly 
matches the prediction given by the conformal field theory of the orbifold. 
\end{reptheorem}

%We also see that this is a minimum number of deformations, not just minimum for those corresponding to case 1).  This is because every other triangulation 
%we could get comes from a sequence of flops.  But the only curves we can flop are the (-1,-1)-curves, and any such flop cannot lower the number of 
%deformations.  It ``adds'' 4 deformations, one along each edge bordering the flop going toward the new endpoints.  If there are already arrows going the 
%opposite direction, it subtracts 1 from each of those.  Because of the construction used, there are at most two edges that do this.  Therefore number of 
%deformations either increases or stays the same-> minimum.

\section{Future Directions}

We would like to understand more about the resolutions where the count does not agree.  This can happen in two ways.  First, there are resolutions where 
each singlet does not correspond to a deformation.  This can be seen from the proof of proposition 7.5, where we saw that the non-negativity conditions 
were not necessarily satisfied.  Second, there are resolutions which have deformations that do not satisfy the conditions for the singlets.  For example, 
we can look at the resolution of $\CC^3/\ZZ_{11}$ corresponding to the following triangulation.
 \begin{center}
\includegraphics[height=2in]{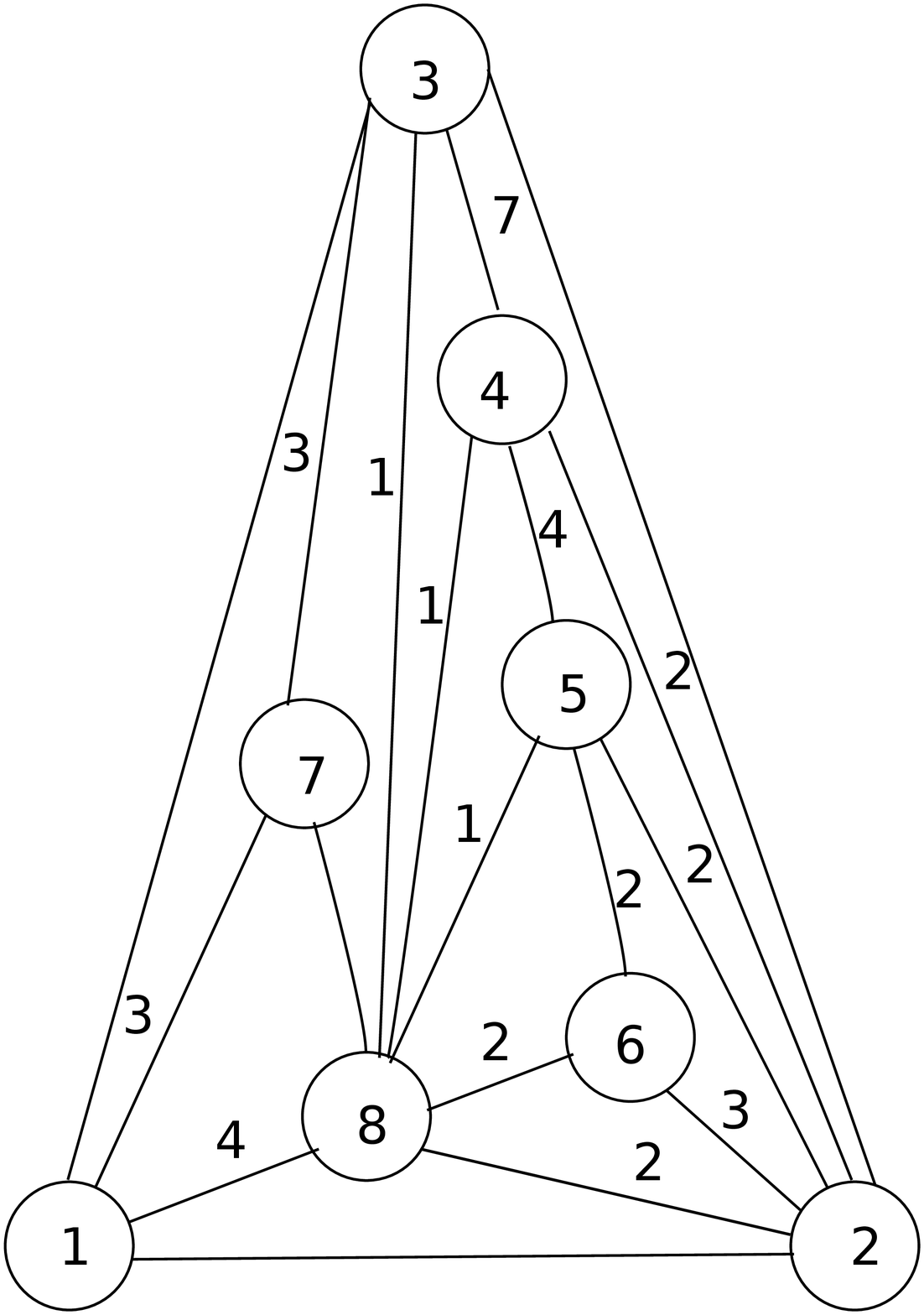}\\
Figure 12.1: Triangulation for Alternate Resolution of $\CC^3/\ZZ_{11}$.\\
\end{center}
We see that this triangulation does not have any deformations coming from $\ext^1(\cO_{D_4}(q_4),\cO_{D_7}(q_7))$, so there is a singlet that does not 
correspond to any deformation.  It also has a deformation from $\ext^1(\cO_{D_4}(q_4),\cO_{D_8}(q_8))$, which does not come from any singlet, so it 
embodies both of the issues discussed above.  

\section{Acknowledgments}

I would like to thank P.~Aspinwall for the suggestion of the problem, and many useful discussions and comments.  I also thank N.~Addington 
for his comments on the organization of the paper.  This work was partially
supported by NSF grant DMS--0905923.  Any opinions, findings, and
conclusions or recommendations expressed in this material are those of
the author and do not necessarily reflect the views of the National
Science Foundation.

\bibliography{biblio}

\end{document}